\documentclass[12pt]{amsart}

\usepackage{bbm}
\usepackage{mathrsfs}
\usepackage{amsfonts}
\usepackage{amsmath}
\usepackage{amssymb}

\usepackage{graphicx}

\usepackage{latexsym}
\usepackage{cases}
\usepackage{amscd}
\usepackage{amsfonts}

\usepackage[all]{xy}
\usepackage{xy}

\usepackage{amsfonts,amsmath, amssymb}

\newtheorem{theorem}{Theorem}[section]

\newtheorem{lemma}[theorem]{Lemma}
\newtheorem{proposition}[theorem]{Proposition}
\newtheorem{remark}[theorem]{Remark}
\newtheorem{definition}[theorem]{Definition}
\newtheorem{example}[theorem]{Example}

\begin{document}

\title{Infinite Dimensional Manifolds from a New Point of View}

\author {Xianzu Lin }

\date{ }
\maketitle
   {\small \it College of Mathematics and Computer Science, Fujian Normal University, }\\
    \   {\small \it Fuzhou, {\rm 350108}, China;}\\
     {\small \it Institute of Mathematics, Academy of Mathematics and Systems
   Science, }\\
    \   {\small \it Beijing {\rm 100190}, China}\\
      \              {\small \it Email: linxianzu@126.com}
\begin{abstract}

In this paper we propose a new treatment about infinite dimensional
manifolds, using the language of categories and functors. Our
definition of infinite dimensional manifolds is a natural
generalization of finite dimensional manifolds in the sense that de
Rham cohomology and singular cohomology can be naturally defined and
the basic properties (Functorial Property, Homotopy Invariant,
Mayer-Vietoris Sequence) are preserved. In this setting we define
the classifying space $BG$ of a Lie group $G$ as an infinite
dimensional manifold. Using simplicial homotopy theory and the
Chern-Weil theory for principal $G$-bundles we show that de Rham's
theorem holds for $BG$. Finally we get, as an unexpected byproduct, two new simplicial set models for the
classifying spaces of compact Lie groups; they are totally different from the classical models constructed by Milnor, Milgram, Segal and Steenrod.

\


\

2000 MR Subject Classification: 55U10,18A05
\end{abstract}

\section{Introduction}

Usually, infinite dimensional manifolds are defined to be
paracompact spaces $X$ modelled on some topological vector space
$E$; $E$ may be a Fr\'{e}che space, a Banach space or a Hilbert
space. More explicitly, $X$ is a paracompact space covered by an
atlas of open subsets $\{U_{\alpha}\}$ each of which is homeomorphic
to an open set $E_{\alpha}$ of $E$ by a given homeomorphism
$\phi_{\alpha}:U_{\alpha}\rightarrow E_{\alpha}$. The transitive
functions between charts are assumed to be infinitely differentiable. The meaning of infinitely differentiable is too
complicated to be given here; we refer the reader to \cite{h,m}for
details. In this setting differential forms can be defined and de
Rham's theorem holds under some mild conditions.

In this paper we propose a new treatment about infinite dimensional
manifolds. Our treatment differs from the usual one in three
respects. First, we adopt the language of categories and functors and
our definition of infinite dimensional manifolds is a natural
generalization of finite dimensional manifolds in the sense that de
Rham cohomology and singular cohomology can be naturally defined and
the basic properties (Functorial Property, Homotopy Invariant,
Mayer-Vietoris Sequence) are preserved. Second, we need no topology
except the topology of finite dimensional manifolds. Thus we would
not talk about the topology of an infinite dimensional manifold,
instead we still have the homotopy type of an infinite dimensional
manifold. Finally we can naturally define many classes of infinite dimensional
manifolds including space of smooth mappings between finite
dimensional manifolds, group of diffeomorphisms of a finite
dimensional manifold, space of connections on a principal G-bundle
and classifying space of Lie group.

In \S2 we give the definition of infinite dimensional manifolds and
several examples. \S3 is concerned with the cohomology theory of infinite
dimensional manifolds. In \S4 we prove that de Rham's
theorem holds for
the classifying space of a compact Lie group $G$. During the proof we get, as an unexpected byproduct, two new simplicial set models for the
classifying spaces of compact Lie groups; they are totally different from the classical models constructed by Milnor, Milgram, Segal and Steenrod (cf.\cite{mi,m1,se,st}). Throughout this paper, we assume that the reader is familiar with some standard results in simplicial homotopy theory. The standard references are \cite{c,gj}.

\section{Definitions and Examples}
In this section we give the definition of infinite dimensional
manifolds. We assume the reader is familiar with the elementary
language of category theory for which \cite{mac} is a good reference.

Let $\mathbb{R}^{n}_{k}\subset\mathbb{R}^{n}$ be the subspace
$$\mathbb{R}^{n}_{k}=\{(x_{1},\cdots x_{n})\in\mathbb{R}^{n}\mid x_{i}\geq0,\ 1\leq i\leq k\}$$
A function $f$ on an open subset $U\subset\mathbb{R}^{n}_{k}$ is
$smooth$ if and only if $f$ can be extended to a $C^{\infty}$ function on an
open subset of $\mathbb{R}^{n}$. A map $g$ from an open subset
$U\subset\mathbb{R}^{n}_{k}$ to an open subset
$U'\subset\mathbb{R}^{m}_{l}$ is $smooth$ if $f_{i}\circ g$ is
smooth for each coordinate function $f_{i}$, $1\leq i\leq m$.
\begin{definition}
A chart (with corners) on a topological space $X$ is a homemorphism
$\phi:U\rightarrow X_{\alpha}$ of an open subset $U$ of
$\mathbb{R}^{n}_{k}$ with an open subset $X_{\alpha}$ of $X$. Two
charts $\phi_{1}$, $\phi_{2}$ are said to be $compatible$ if the
functions $\phi_{2}^{-1}\phi_{1}$ and $\phi_{1}^{-1}\phi_{2}$ are smooth.
\end{definition}
\begin{definition}
A $smooth \ atlas$ (with corners) on a space $X$ is an countable
open cover $\{X_{\alpha}\}$ of $X$ together with a collection of compatible
charts $\phi_{\alpha}:U_{\alpha}\rightarrow X_{\alpha}$. A $smooth\
manifold \ with \ corners$ is a paracompact Hausdorff space $X$
together with a maximal smooth atlas, i.e., a smooth structure.
\end{definition}
As usual we can define $smooth$ $map$ and $smooth$ $homeomorphism$ between smooth manifolds with corner.
We also have the concepts of closed submanifold and open submanifold, and the Whitney embedding theorem still holds in this setting; the proof is routine. From now on, we assume that all smooth manifolds are Hausdorff, paracompact and have a countable base.

\begin{definition}
Let $\mathscr{M}$ the category of smooth manifolds with corner and
smooth maps and let $(Set)$ be the category of sets. A smooth
functor is a contravariant functor
$$F:\mathscr{M}^{op}\rightarrow(Set).$$
A smooth functor $F$ is called separated (resp. an infinite
dimensional manifold) if it satisfies the following condition: Assume that
$\bigcup_{i}V_{i}$ is an open covering of $X$ in $\mathscr{M}$, and
for each $i$ there is an $\alpha_{i}\in F(V_{i})$ such that
$\alpha_{i}|V_{i}\cap V_{j}=\alpha_{j}|V_{i}\cap V_{j}$ for
each $i,j$, then there is at most (resp. exactly) one $\alpha\in
F(X)$ such that $\alpha|V_{i}=\alpha_{i}$.
\end{definition}

\begin{remark}
In the language of sheaf theory (cf. \cite{ar,ta}), a smooth functor is just a presheaf on the site of smooth manifolds, and a separated smooth functor (resp. an infinite dimensional manifold) is a separated presheaf (resp. a sheaf).
\end{remark}

\begin{definition}
Let $F$, $G$ be smooth functors, a $smooth\ map$ from $F$ to $G$ is
a natural transformation from $F$ to $G$. Denote by
$\widetilde{\mathscr{M}}$ the category of smooth functors and smooth
maps, and denote by $\widetilde{\mathscr{M}}'$ and
$\mathbb{M}$ the full subcategory of separated smooth
functors and infinite dimensional manifolds respectively.
\end{definition}
\begin{lemma}\label{lm:a}
The inclusion functor
$i:\widetilde{\mathscr{M}}'\hookrightarrow\widetilde{\mathscr{M}}$ has a
left adjoint functor
$\varsigma:\widetilde{\mathscr{M}}\rightarrow\widetilde{\mathscr{M}}'$,
and the inclusion functor
$i':\mathbb{M}\hookrightarrow\widetilde{\mathscr{M}}'$
has a left adjoint functor
$\varsigma':\widetilde{\mathscr{M}}'\rightarrow\mathbb{M}$
\end{lemma}

This lemma is standard in sheaf theory (cf.\cite{ta}), we give a proof for later use.
\begin{proof}
Given a smooth functor $F$, we construct $\varsigma(F)$ as follows.
Let $M$ be an object in $\mathscr{M}$, we say that $\alpha,\beta\in
F(M)$ are equivalent if there is an open covering $\bigcup_{i}M_{i}$
of $M$ such that $\alpha|M_{i}=\beta|M_{i}$ for each $i$.
Define $\varsigma(F)(M)$ to be the equivalent classes of elements of
$F(M)$. It is clear that $\varsigma(F)(M)$ is functorial in $M$ and
$\varsigma(F)$ is a separated smooth functor.

Given a separated smooth functor $G$, we construct $\varsigma'(G)$
as follows. Let $M$ be an object in $\mathscr{M}$, a $local$ $datum$ of
$G$ on $M$ consists of an open covering $\bigcup_{i}M_{i}$ of $M$,
and for each $i$ an $\alpha_{i}\in G(M_{i})$, satisfying
$\alpha_{i}|M_{i}\cap M_{j}=\alpha_{j}|M_{i}\cap M_{j}$ for
any $i$, $j$. Two local data $(\bigcup_{i}M_{i}, \alpha_{i})$ and
$(\bigcup_{j}N_{j}, \beta_{j})$ of $M$ on $G$ are said to be
$equivalent$ if $\alpha_{i}|M_{i}\cap N_{j}=\beta_{j}|
M_{i}\cap N_{j}$ for any $i$, $j$. Define $\varsigma'(G)(M)$ to be
the set of equivalent classes of local data of $G$ on $M$. Then
$\varsigma'(G)(M)$ is functorial in $M$ and $\varsigma'(G)$
is an infinite dimensional manifold. One checks that $\varsigma$ and $\varsigma'$ are adjoint functors of
$i$ and $i'$ respectively.
\end{proof}

\begin{definition}
Let $F$  be an infinite dimensional manifold, a $submanifold$
of $F$ is an infinite dimensional manifold $G$ such that, for
each $M$ in $\mathscr{M}$, $G(M)$ is a subset of $F(M)$.
\end{definition}

Given two submanifolds $G$, $H$ of an infinite dimensional
manifold $F$, define $G \sqcup H$ in $\widetilde{\mathscr{M}}'$ by
$G\sqcup H(M)=G(M)\cup H(M)$. Set $G\cup H=\varsigma'(G \sqcup H)$,
then $G\cup H$ is also a submanifold of $F$, called the $union$ of
$G$ and $H$. The $intersection$ of $G$ and $H$ is the submanifold
$G\cap H$ defined by $G\cap H(M)=G(M)\cap H(M)$. The $product$ of two infinite dimensional
manifolds $F$, $G$ is the infinite dimensional
manifold $F\times G$ defined by $F\times G(M)=F(M)\times G(M)$.

To each $M$ in $\mathscr{M}$, we can associate a contravariant
functor
$$h_{M}:\mathscr{M}^{op}\rightarrow(Set)$$
which sends $N$ in $\mathscr{M}$ to the set
$Hom_{\mathscr{M}}(N,M)$; if $\alpha:N\rightarrow N^{'}$
is an arrow in $\mathscr{M}$, then $h_{M}(\alpha)$ is defined by composition with $\alpha$. We say that $h_{M}$ is a cofunctor
represented by $M$. It is clear that $h_{M}$ is an infinite
dimensional manifold in our sense. An arrow $f:M\rightarrow N$ in
$\mathscr{M}$ yields a smooth map, i.e., a natural transformation
$h_{f}:h_{M}\rightarrow h_{N}$.
In this setting, the Yoneda Lemma states that,
\begin{lemma}
Let $M$ and $N$ be objects in $\mathscr{M}$, then the function
$$Hom_{\mathscr{M}}(M,N)\rightarrow Hom_{\mathbb{M}}(h_{M},h_{N})$$
that sends $f:M\rightarrow N$ to $h_{f}$ is bijective.
\end{lemma}
From now on we always identify each $M$ in $\mathscr{M}$ with the
infinite dimensional manifold $h_{M}$ in $\widetilde{\mathscr{M}}$.
Thus $\mathscr{M}$ is a full subcategory of
$\mathbb{M}$.

\begin{example}
To any $M$, $N$ in $\mathscr{M}$, we can associate an infinite
dimensional manifold $C^{\infty}(M,N)$ by sending $L\in\mathscr{M}$
to $Hom_{\mathscr{M}}(L\times M,N)$, if $f:L\rightarrow L'$ is an
arrow in $\mathscr{M}$, $C^{\infty}(M,N)(f)$ is defined by
composition with $f\times\mathbb{I}:L\times M\rightarrow L'\times
M$.
\end{example}

\begin{example}
For any $M$ in $\mathscr{M}$, define $Diff(M)\in
\mathbb{M}$ by sending $N\in\mathscr{M}$ to the set
$$\{f\in Hom_{\mathscr{M}}(N\times M,M)\mid f|\{x\}\times M \ is \ a \ diffeomorphism \ for \ each \ x\in N
\}.$$  For an arrow $f:L\rightarrow L'$ in $\mathscr{M}$,
$Diff(M)(f)$ can be defined analogously.
\end{example}

\begin{example}
Let $G\times M\rightarrow M$ be a smooth action of a Lie group $G$ on $M\in\mathscr{M}$. Then the group $Hom_{\mathscr{M}}(N,G)$ acts on $Hom_{\mathscr{M}}(N,M)$. Define $M^G\in
\mathbb{M}$ by sending $N\in\mathscr{M}$ to the quotient set $Hom_{\mathscr{M}}(N,M)/Hom_{\mathscr{M}}(N,G)$.
\end{example}

\begin{example}
Let $\xi=\{E\rightarrow M\}$ in $\mathscr{M}$ be a vector bundle or
a principal $G$-bundle, define $Sec(\xi)\in
\mathbb{M}$ by sending each $N$ in
$\mathscr{M}$ to the set of smooth section of
$P^{*}\xi=\{P^{*}E\rightarrow M\times N\}$ where $P:M\times
N\rightarrow M$ is the projection. For each arrow $f:N\rightarrow
N'$, $Sec(\xi)(f)$ is defined by the pullback of sections along
$f\times\mathbb{I}:M\times N\rightarrow M\times N'$.
\end{example}

\begin{example}
Let $G$ be a compact Lie group. A $G$-fiber in $\mathbb{R}^{n}$ is
a closed sub-manifold $M$ of $\mathbb{R}^{n}$ with a smooth
$G$-action such that for any $x\in M$, the map $f_{x}:G\rightarrow
M$ defined by $f_{x}(g)=g\cdot x$ is a diffeomorphism.  Let $G_n$
be the set of all $G$-fibers in $\mathbb{R}^{n}$. Then we have
$G_n\subset G_{n+1}\subset G_{n+2}\subset \cdots$ Set
$$G_\infty=\bigcup_n G_n.$$ For any $M\in \mathscr{M}$, a map $f:M\rightarrow G_\infty$
is said to be $smooth$ if

 $({\bf 1})$ $f(M)\subset G_n$ for some $n$;

 $({\bf 2})$  $P_f:=\bigcup_x\{x\}\times
f(x)$ is a closed submanifold of $M\times\mathbb{R}^{n}$;

 $({\bf
3})$ the projection $P_f\rightarrow M$ is a principal $G$-bundle.

Now the classifying space $BG\in\mathbb{M}$   of $G$ is defined by
sending each $M$ in $\mathscr{M}$ to the set of smooth maps
$f:M\rightarrow G_\infty$ with a connection on the principal
$G$-bundle $P_f$; for an arrow $f:M\rightarrow N$ in
$\mathscr{M}$, $BG(f)$ is defined by composition with $f$. It is
easy to see that $BG$ is indeed an infinite dimensional manifold.
Our definition of $BG$ is suggested by the existence of universal
connection on the classifying space of Lie group
(cf.\cite{nr,nr1}).
\end{example}

\begin{example}
Let $\xi=\{E\rightarrow M\}$ in $\mathscr{M}$ be a principal
$G$-bundle, the $space\ of\ connections$ of $\xi$ is an infinite
dimensional manifold $Con(\xi)$ defined as
follows.
For each object $N$ in $\mathscr{M}$, $Con(\xi)(N)$ is the set of
connections on the principal $G$-bundle $P^{*}(\xi)$, where
$P:M\times N\rightarrow M$ is the projection. For an arrow
$f:N\rightarrow N'$ in $\mathscr{M}$, $Con(\xi)(f)$ is defined by
the pullback of connections along $f\times\mathbb{I}:M\times
N\rightarrow M\times N'$.
\end{example}

\begin{example}
For each $n\geq0$, $\mathcal {A}^{n}\in
\mathbb{M}$ is defined by sending each $M$ in
$\mathscr{M}$ to the set of differential $n$-forms on $M$; for each
arrow $f:M\rightarrow N$ in $\mathscr{M}$, $\mathcal {A}^{n}(f)$ is
the induced map $f^{*}:\mathcal {A}^{n}(N)\rightarrow\mathcal
{A}^{n}(M).$ The exterior derivative induces a smooth map
$d^{n}:\mathcal {A}^{n}\rightarrow\mathcal {A}^{n+1}$, and the
exterior product induces a smooth map $\wedge:\mathcal
{A}^{m}\times\mathcal {A}^{n}\rightarrow\mathcal {A}^{m+n}$.
\end{example}


\section{Cohomology Theory}
Let $F$ be a smooth functor, a differential $ n$-form on $F$ is
a smooth map from $F$ to $\mathcal {A}^{n}$. Denote by $\mathcal
{A}^{n}(F)$ the real vector space of differential n-forms on $F$,
the smooth map $d^{n}:\mathcal {A}^{n}\rightarrow\mathcal {A}^{n+1}$ induces a differential operator
$d^{n}(F):\mathcal {A}^{n}(F)\rightarrow\mathcal {A}^{n+1}(F)$.
The complex $\mathcal
{A}^{*}(F)$ together with the differential operator $d^{*}(F)$
is called the $de\ Rham\ complex$ on $F$.
Applying the smooth map $\wedge:\mathcal {A}^{m}\times\mathcal
{A}^{n}\rightarrow\mathcal {A}^{m+n}$ we see that the
$\mathcal {A}^{*}(F)$ has a natural structure of commutative
differential graded algebra, the corresponding cohomology ring
$H^{*}_{de}(F)$ is called the $de\ Rham\ cohomology$ of $F$.
\begin{definition}
Let $F$ be a smooth functor, the $singular\ complex$ $S(F)$ is the
simplicial set given by
$$ S(F)_{n}=F(|\triangle^{n}|)=Hom_{\widetilde{\mathscr{M}}}(|\triangle^{n}|,F),$$
where
$$|\triangle^{n}|=\{(t_{0},\cdots,t_{n})\in\mathbb{R}^{n+1}\mid\sum_{i=0}^{n}t_{i}=1,\ t_{i}\geq0 \}$$
is a manifold with corner. The $singular\ cohomology\ ring$
$H^{*}_{si}(F; R)$ of $F$ with coefficient ring $R$  is the
cohomology ring $H^{*}(S(F); R)$; we can analogously define singular
(co)chain complex and the homology (homotopy) groups of $F$.
\end{definition}
Thus, although we have no $topology$ of an infinite
dimensional manifold, we can still talk about the homotopy type of
it. In the following three propositions $H^{*}(F)$ will be
$H^{*}_{de}(F)$ or $H^{*}_{si}(F; R)$.
\begin{proposition}
({\bf Functorial Property}) Each smooth map $f:F\rightarrow G$ in
$\widetilde{\mathscr{M}}$ induces a ring homomorphism
$f^{*}:H^{*}(G)\rightarrow H^{*}(F)$. If $g:G\rightarrow L$ is
another smooth map in $\widetilde{\mathscr{M}}$, then
$(gf)^{*}=f^{*}g^{*}$.
\end{proposition}
Let $\mathbb{I}\in\mathscr{M}$ be the unit interval $[0,1]$ and let
$f,g:F\rightarrow G$ be two smooth maps in $\widetilde{\mathscr{M}}$.
We say that $f$ is $smooth\ homotopic$ to $g$ if there is a smooth
map $h:F\times\mathbb{I}\rightarrow G$ such that $h(\cdot,0)=f$ and
$h(\cdot,1)=g$.
\begin{proposition}
({\bf Homotopy Invariant})Let $f,g:F\rightarrow G$ be two smooth
maps in $\widetilde{\mathscr{M}}$, if $f$ is smooth homotopic to
$g$, then $f^{*}=g^{*}:H^{*}(G)\rightarrow H^{*}(F)$.
\end{proposition}
\begin{proof}
The proof is routine, we give proof only in the case $H^{*}=H^{*}_{de}$. Let
$h:F\times\mathbb{I}\rightarrow G$ be a smooth map such that
$h(\cdot,0)=f$ and $h(\cdot,1)=g$. It suffices to construct a chain
homotopy between $f^{*},g^{*}:\mathcal {A}(F)\rightarrow\mathcal
{A}(G)$. Given $M$ in $\mathscr{M}$, any $n$-form $w$ on
$M\times\mathbb{I}$ can be uniquely written as $$w=w'(t)\wedge dt+w''(t)$$
where $w'(t)$ (resp.$w''(t)$) is an $n-1$-form (resp. $n$-form) on $M$ (smoothly dependent on $t$). By
the uniqueness $w'$ and $w''$ are functorial in $M$. Thus
for any $w\in\mathcal {A}^{n}(G)$, we can set $h^{*}(w)=w'(t)\wedge dt+w''(t)$ as above, where $w'(t)$ (resp.$w''(t)$) is an $n-1$-form (resp. $n$-form) on $G$. Define $s(w)=\int_{0}^{1}w'(t)dt$, then we have
$ds+sd=g^{*}-f^{*}$ as in the finite dimensional case. Thus $s$ is a
chain homotopy between $f^{*},g^{*}$ and the proof of the
proposition is done.
\end{proof}

\begin{proposition}For any $F$ in $\widetilde{\mathscr{M}}$ and $G$ in $\widetilde{\mathscr{M}}'$,
the canonical maps $h:F\rightarrow\varsigma(F)$ and
$h':G\rightarrow\varsigma'(G)$ induce isomorphisms on singular cohomology and de Rham cohomology.
\end{proposition}
\begin{proof}We give the proof only in the case of singular cohomology.
First we show that $h$ and $h'$ induce isomorphisms on singular homology. Let
$C_{*}(F)$ be the singular chain complex of $F$ and let $S$ be the
subdivision operator $S(F):C_{*}(F)\rightarrow C_{*}(F)$, it was
shown in \cite[p.121]{ha} that $S(F)$ is chain homotopic to the
identity map. Let $H'_{*}(F)$ be the homology groups of the complex $ker\ h_{*}$, where $ker\ h_{*}$ is the kernel of the surjective
homomorphism $h_{*}:C_{*}(F)\rightarrow C_{*}(\varsigma(F))$. As $S$
induces isomorphisms on $H_{*}(F)$ and $H_{*}(\epsilon(F))$, it also
induces isomorphisms on $H'_{*}(F)$. On the other hand, by the definition of $\varsigma(F)$ in the proof of Lemma \ref{lm:a},  we see that for any
singular chain $\alpha\in ker\ h_{*}$ we have $S^{n}(\alpha)=0$ for
sufficiently large $n$. Hence
we have $H'_{*}(F)\cong0$ and $h$ induces isomorphism on singular
homology.

Now we show that for any separated smooth functor $G$ the canonical map
$h':G\rightarrow\varsigma'(G)$ induces isomorphism on singular
homology. We only need to repeat the argument of the previous paragraph
except that the kernel of $h_{*}$ should be replaced by the
cokernel of $h'_{*}$, the details are omitted. By the universal
coefficients theorem the same is true for cohomology.
\end{proof}

\begin{proposition}
({\bf Mayer-Vietoris Sequence})Let $F$  be an infinite
dimensional manifold and let $G$, $H$ be two submanifolds of $F$
such that $G\cup H=F$, then we have the following long exact
sequence
$$\xymatrix@C=0.5cm{
 \cdots  \ar[r] & H^{n-1}_{si}(G\cap H; R) \ar[r]  & H^{n}_{si}(F; R)\ar[r]
 & H^{n}_{si}(G; R)\oplus H^{n}_{si}(H; R)\\
\ar[r] & H^{n}_{si}(G\cap H; R) \ar[r]  & H^{n+1}_{si}(F; R)\ar[r] & \cdots }$$
\end{proposition}
\begin{proof}
From the short exact sequence of chain complexes formed by
$$\xymatrix@C=0.5cm{
0 \ar[r] & C^{n}_{si}(G \sqcup H; R) \ar[r]  & C^{n}_{si}(G; R)\oplus
C^{n}_{si}(H; R)
 \ar[r] &C^{n}_{si}(G\cap H; R)\ar[r] &0
 }$$
 we have the corresponding long exact sequence of cohomology
$$\xymatrix@C=0.5cm{
 \cdots  \ar[r] & H^{n-1}_{si}(G\cap H; R) \ar[r]  & H^{n}_{si}(G \sqcup H; R)
 \ar[r] & H^{n}_{si}(G; R)\oplus H^{n}_{si}(H; R)\\
 \ar[r]& H^{n}_{si}(G\cap H; R)\ar[r]  & H^{n+1}_{si}(G \sqcup H; R)\ar[r] &\cdots}.$$
 By the previous proposition the inclusion $G \sqcup H\hookrightarrow G\cup H=F$ induces isomorphism on
 cohomology, thus the proof is done.
\end{proof}

\begin{remark}
For any topological space $X$, we can define an infinite
dimensional manifold $h_{X}$ as follows. For each object $N$ in
$\mathscr{M}$, $h_{X}(N)$ is the set of all continue maps from $N$
to $X$. For an arrow $f:N\rightarrow N'$ in $\mathscr{M}$,
$h_{X}(f)$ is defined by composition with $f$. Let $A\cup B$ be an
open covering of $X$, then $h_{A}$ and $h_{B}$ are submanifolds of
$h_{X}$ with $h_{A}\cup h_{B}=h_{X}$. Applying the above
proposition to $(h_{X},h_{A}, h_{B})$, we get the usual Mayer-Vietoris Sequence for the triple $(X;
A,B)$.
\end{remark}

Let $F$ be an infinite dimensional manifold. To each differential
n-form $\omega:F\rightarrow\mathcal {A}^{n}$ on $F$, we will associate a singular
n-cochain $\widetilde{\omega}$ on $F$ as follow. For each $\alpha\in
F(\triangle^{n})$, $\omega(\alpha)$ is a
 differential n-form (in the usual sense) on $\triangle^{n}$; set
 $\widetilde{\omega}(\alpha)=\int_{\triangle^{n}}\omega(\alpha)$.

It is clear that this assignment $\omega\rightarrow
\widetilde{\omega}$ induces a homomorphism from the de Rham complex
to the singular cochain complex with real coefficient, denote by
$\mathfrak{R}(F):H^{*}_{de}(F)\rightarrow H^{*}_{si}(F; \mathbb{R})$
the induced homomorphism of cohomology. The famous de Rham theorem states
that $\mathfrak{R}(M) $ is an isomorphism for each $M\in\mathcal {M}$.
In the general case we don't know under what conditions will $\mathfrak{R}(F)$ be an
isomorphism. Our definition of infinite dimensional manifolds is
incomplete until we find the suitable conditions. We refer the
reader to \cite{m3} for the proof of de Rham theorem of infinite
dimensional manifolds in the usual sense.
\begin{remark}
Note that in this section all the definitions, results and their proofs coincide
with the usual ones in the finite dimensional case except the
definition of singular (co)homology. When
$F\in\widetilde{\mathscr{M}}$ is represented by $M$ in
$\mathscr{M}$, $H^{*}_{si}(F; \mathbb{R})$ is not the singular
cohomology of $M$ but a cohomology defined by
smooth singular (co)chains. The equivalence of these two
cohomologies has been proved in \cite{ei}.
\end{remark}

\section{Classifying Spaces of Lie Groups}
In this section we show that the de Rham theorem is valid for
$BG\in\mathbb{M}$ where $G$ is a compact Lie group in
$\mathscr{M}$.

Consider $BG'\in\widetilde{\mathscr{M}}$ which sends each $M$ in
$\mathscr{M}$ t the set of smooth maps $f:M\rightarrow G_\infty$
. The map $p':BG\rightarrow BG'$ is defined by neglecting
the connections.

\begin{theorem}\label{de}
The geometric realization $|S(BG)|$ of the simplical set $S(BG)$ is a classifying space of $G$.
\end{theorem}

\begin{proof}
We divide the proof into four lemmas.

\begin{lemma}
$S(BG')$ is a connected fibrant.
\end{lemma}
\begin{proof}
It is easy to see that $S(BG')$ is a connected simplicial set.
Thus it remains to show that given a commutative diagram in the
category of simplicial set $\mathbb{S}$:
\begin{equation}
\xymatrix{
  \Lambda^{n}_{k} \ar[d]_{i} \ar[r]
                & S(BG') \ar[d]^{p}  \\
  \triangle^{n} \ar@{.>}[ur] \ar[r]
                & \ast    ,      }
\end{equation}

The standard $n$-simplicial set $\triangle^{n}$ is the contravariant functor represented by $\mathbbm{n}$, i.e.,
$$\triangle^{n}=hom_{\Delta}(\cdot,\mathbbm{n}).$$Write $\iota_{n}=\mathbbm{1}_{n}\in
hom_{\Delta}(\mathbbm{n},\mathbbm{n})=\triangle^{n}_{n}$.
The $k$-th
horn $\Lambda^{n}_{k}\subset\triangle^{n}$

   there is a map $\theta:\triangle^{n}\rightarrow S(BG')$
   (the dotted arrow) making the diagram commute.
   Equivalently given a principal $G$-bundle (in $\mathscr{M}$) on $|\Lambda^{n}_{k}|$
   we want to extend it to a principal $G$-bundle (in $\mathscr{M}$) on
   $|\triangle^{n}|$. As the inclusion $|\Lambda^{n}_{k}|\hookrightarrow|\triangle^{n}|$ is a homotopy
   equivalence, such an extension always exists.
\end{proof}

\begin{lemma}
$|S(BG')|$ is a classifying space of $G$.
\end{lemma}

\begin{proof}
Let $\mathbb{BG}$ (a $CW$ complex) be a classifying space of $G$.
For any finite subcomplex $K$ of $|S(BG')|$, there is a canonical
principal $G$-bundle $\xi_{K}$ over $K$ with projection
$E_{K}\rightarrow K$ (glue together principal $G$-bundles on each
simplice), hence a classifying map $l_{K}:K\rightarrow\mathbb{BG}$.
As $\{l_{K}\}$ are compatible (up to homotopy) under
inclusions of finite subcomplexes of $|S(BG')|$, by the homotopy extension property of CW complexes (cf.\cite{ha}), they induce a map
$l:|S(BG')|\rightarrow \mathbb{BG}$ such that the restriction of $l$ to $K$ is homotopic to $l_{K}$ for each finite subcomplex $K$. In
order to prove this lemma, it suffices to show that $l$ induces
isomorphism $l_{*}$ on homotopy groups.

First we show that $l_{*}$ is surjective. Fix a faithful
representation $i:G\rightarrow U(n)$ (for the existence of such a representation, see \cite{bt}, p.136), and consider the principal $G$-bundle $\delta_{m}$
 with projection $U(m+n)/U(m)\rightarrow U(m+n)/(U(m)\times G)$ for each $m$.
Choose a triangulation of $U(m+n)/(U(m)\times G)$, then it induces a
classifying map (by the definition of $BG'$) $k:U(m+n)/(U(m)\times
G)\rightarrow |S(BG')|$ such that the composition $l\cdot k$ is the classifying map
for $\delta_{m}$. As $U(m+n)/U(m)$ is $2m$-connected,
the classifying map $l\cdot k$ induces isomorphism on
$\pi_{i}$ for $i<2m$ (cf.\cite{sw}, p.202), hence $l_{*}$ is surjective on homotopy groups.

Now we show that $l_{*}$ is injective. Let
$\alpha\in\pi_{n}(|S(BG')|)$ satisfying $l_{*}(\alpha)=0$, we want
to show that $\alpha=0$. As $S(BG')$ is a fibrant, we can represent
$\alpha$ by a simplicial map $\alpha:\triangle^{n}\rightarrow S(BG')$ such that
$\alpha|\partial\triangle^{n}$ is the constant map $\partial\triangle^{n}\rightarrow
*\hookrightarrow S(BG')$ (where $*$ is represented by the trivial $G$-bundle $G\rightarrow*$).

By the definition of $BG'$, $\alpha$ is represented by a principle
$G$-bundle $\xi$ (in $\mathscr{M}$) on $|\triangle^{n}|$ satisfying $\xi|\partial\triangle^{n}$ is trivial, i.e. a
principle $G$-bundle $\xi$ on
$|\triangle^{n}/\partial\triangle^{n}|\cong S^{n}$. Now
$l_{*}(\alpha)=0$ implies $\xi$ is isomorphic to the trivial
$G$-bundle
$$G\times|\triangle^{n}/\partial\triangle^{n}|\longrightarrow|\triangle^{n}/\partial\triangle^{n}|.$$
Now condition ({\bf 5}) of Theorem \ref{main} implies that there is a principal
$G$-bundle $\delta$ (in $\mathscr{M}$) over
$|\triangle^{n}|\times\mathbb{I}$ together with a trivialization of
$\delta|(|\partial\triangle^{n}|\times\mathbb{I})$, such that
$\delta|(|\triangle^{n}|\times\{0\})=\xi$ and
$\delta|(|\triangle^{n}|\times\{1\})$ is the trivial $G$-bundle.
This yields a homotopy from $\alpha$ to the constant map, thus
$\alpha=0$.
\end{proof}

\begin{lemma}\label{tt}
Let $\mathbb{R}_{+}^{n}=\{(t_{1},\cdots,t_{n})\in\mathbb{R}^{n}|t_{i}\geq0
  \}$
  and $H_{i}=\{(t_{1},\cdots,t_{n})\in\mathbb{R}_{+}^{n}\mid t_{i}=0
  \}.$
  Given a smooth function $f_{i}$ on $H_{i}$ for each $i$ satisfying $f_{i}|H_{i}\cap H_{j}=f_{j}|H_{i}\cap
 H_{j}$, there always exists a smooth function $f$ on $\mathbb{R}_{+}^{n}$ such
 that $f|H_{i}=f_{i}$.
\end{lemma}
\begin{proof}
Extend $f_{n}$ to a smooth function $f'_{n}$ on $\mathbb{R}_{+}^{n}$
by setting
$f'_{n}((t_{1},\cdots,t_{n})=f_{n}((t_{1},\cdots,t_{n-1},0)$. If we
can find a smooth function $f''$ satisfying
$f''|H_{i}=f_{i}-f'_{n}$, then $f''+f'_{n}$ is the desired $f$. Thus
it suffices to prove this claim when $f_{n}\equiv0$. Repeating this
argument we see that it suffices to prove this lemma when
$f_{i}=0$ for any $i\leq n$; in this case the lemma is trivial.
\end{proof}

\begin{lemma}
The induced simplicial map $p:S(BG)\rightarrow S(BG')$ is a
fibration.
\end{lemma}
\begin{proof}
It suffices to prove that for every commutative diagram in
$\mathbb{S}$

\begin{equation}
\xymatrix{
  \Lambda^{n}_{k} \ar[d]_{i} \ar[r]
                & S(BG) \ar[d]^{p}  \\
  \triangle^{n} \ar@{.>}[ur] \ar[r]
                & S(BG')           }
\end{equation}

   there is a map $\theta:\triangle^{n}\rightarrow S(BG)$
   (the dotted arrow) making the diagram commute.

 By the definition of $BG$ and $BG'$, the above diagram gives a
 trivial principle $G$-bundle $\xi$ on $|\triangle^{n}|$ and for
 each $i$
 a connection $w_{i}$ on $(d^{i})^{*}(\xi)$ except when $i=k$, satisfying
 $w_{i}=w_{j}$ when restrict to intersection of two faces. To find
a map $\theta:\triangle^{n}\rightarrow S(BG)$ making the diagram
commute is equivalent to find a connection $w$ on $\xi$ such that
when restricted to each $(d^{i})^{*}(\xi)$ ($i\neq k$), $w=w_{i}$.
But in a trivial G-bundle, a connection is just a
 $\mathbbm{g}$-value differential 1-form on the base
  space, where $\mathbbm{g}$ is the Lie algebra of $G$. By a
  linear coordinate transformation and Lemma \ref{tt} we see that there always exists such a connection on  $\xi$
\end{proof}

\begin{lemma}
The fibre $f^{-1}(*)$ is contractible.
\end{lemma}
\begin{proof}
By Lemma 5.1. in \cite[p.190]{gj}, it suffices to show that $f^{-1}(*)$ is
connected and has an extra degeneracy. The $n$-simplices of
$f^{-1}(*)$ is the set of connections on the trivial principal
$G$-bundle $p:G\times|\triangle^{n}|\rightarrow|\triangle^{n}|$,
i.e., the set of
 $\mathbbm{g}$-value differential 1-forms on $|\triangle^{n}|$. Thus there is only
one $0$-simplices {*}. To each $\mathbbm{g}$-value differential
1-form $w$ on $|\triangle^{n}|$ we assign a $\mathbbm{g}$-value
differential 1-form $s_{-1}(w)$ on $|\triangle^{n+1}|$ as follows.

Set
$$\triangle^{+}=\{(t_{0},\cdots,t_{n+1})\in|\triangle^{n+1}|\mid\
t_{0}\leq\frac{1}{2} \}$$ and
$$\triangle^{-}=\{(t_{0},\cdots,t_{n+1})\in|\triangle^{n+1}|\mid\
t_{0}\geq\frac{1}{2} \}.$$
Then $|\triangle^{n+1}|=\triangle^{+}\cup\triangle^{-}$. Define a
projective $f:\triangle^{+}\rightarrow |\triangle^{n}|$ by
$f(t_{0},\cdots,t_{n+1})=(\frac{t_{1}}{1-t_{0}},\cdots,\frac{t_{n+1}}{1-t_{0}})$
and set
$s_{-1}(w)|\triangle^{+}=(e^{4}e^{-(\frac{1}{2}-t_{0})^{-2}})f^{*}(w)$,
$s_{-1}(w)|\triangle^{-}=0$. One checks that $s_{-1}(w)$ is a
$\mathbbm{g}$-value differential 1-form on $|\triangle^{n+1}|$ and
this assignment gives an extra degeneracy on $f^{-1}(*)$. Thus the
proof is done.
\end{proof}
Applying Lemma 4.3, Lemma 4.5, Lemma 4.6 and long exact sequence of homotopy groups in \cite[p.28]{gj} we see that $p':S(BG)\rightarrow S(BG')$ induces isomorphisms on homotopy groups, hence $|S(BG)|$ is also a classifying space of $G$. This completes the proof of Theorem \ref{de}
\end{proof}

\begin{proposition}
$\mathfrak{R}(BG):H^{*}_{de}(BG)\rightarrow H^{*}_{si}(BG; \mathbb{R})$ is an isomorphism. In other words, de Rham's
theorem holds for $BG$.
\end{proposition}
\begin{proof}
By the definition of $BG$, the differential forms on $BG$ are gauge
natural differential forms. Theorem 52.8 in \cite[p.403]{kms}
states that all gauge natural differential forms are the classical
Chern-Weil forms. From Theorem 4.1 and the Chern-Weil theory (we
refer the readers to \cite[Ch.12]{kn} for a detailed exposition about
Chern-Weil theory) we see that
$\mathfrak{R}(BG):H^{*}_{de}(BG)\rightarrow H^{*}_{si}(BG; \mathbb{R})$ is an isomorphism.
\end{proof}
\begin{remark}
We see that $BG$ and $BG'$ are two new (simplicial set) models for the
classifying spaces of compact Lie groups; they are totally different from the classical models constructed by Milnor, Milgram, Segal and Steenrod \cite{mi,m1,se,st}. The referee pointed out that these models have algebraic-geometry analogs for algebraic groups (cf. \cite{lm}).
\end{remark}

\noindent {\bf Acknowledgments}

The author warmly thanks his supervisor Duan Haibao for supports
during the writing of this paper. Thanks are also due to Liu
Xiaobo and Ai Yonghua for their discussions and comments.

\end{document}